\documentclass[11pt]{amsart}

\usepackage[dvipsnames]{xcolor}
\usepackage{amsmath,
amssymb,
amscd,
amsfonts,
array,
comment,
enumitem,
epsfig,
float,
graphicx,
latexsym,
longtable,
mathrsfs,
mathtools,
pdfpages,
psfrag,
rotating,
transparent,
url,
wrapfig,
xr}
\usepackage[linktoc=page]{hyperref}
\hypersetup{
colorlinks,
citecolor=blue,
filecolor=blue,
linkcolor=blue,
urlcolor=blue
}
\usepackage[all]{xy}
\usepackage{amsaddr}
\usepackage[left=1.5in,right=1.5in,top=1.5in,bottom=1.5in+13pt,footskip=0.3in+13pt,headheight=13pt,headsep=0.3in]{geometry}
\usepackage[utf8]{inputenc}
\usepackage{mathtools}
\usepackage{relsize}
\usepackage[hang,flushmargin]{footmisc}
\usepackage[margin=1cm]{caption}

\usepackage{fancyhdr}
\newtheorem{theorem}{Theorem}[section]

\theoremstyle{definition}
\newtheorem{definition}[theorem]{Definition}

\theoremstyle{remark}

\theoremstyle{plain}
\newcommand{\thistheoremname}{}
\newtheorem{genericthm}[theorem]{\thistheoremname}

\newtheorem*{genericthm*}{\thistheoremname}
\newenvironment{namedthm*}[1]
{\renewcommand{\thistheoremname}{#1}%
\begin{genericthm*}}
{\end{genericthm*}}

\newcommand\cF{\mathcal{F}}

\newcommand\cK{\mathcal{K}}

\newcommand\cT{\mathcal{T}}

\newcommand{\sB}{\mathscr{B}}

\newcommand{\sC}{\mathscr{C}}

\newcommand\qu{/\kern-.7ex/} 
\newcommand\lqu{\backslash \kern-.7ex \backslash}

\begin{document}
\sloppy
\raggedbottom
\title{A Polytopal Realization of Higher-Categorical Associahedra}
\author{Spencer Backman, Nathaniel Bottman, and Daria Poliakova}
\begin{abstract}
We describe a polytopal realization of categorical $n$-associahedra.  The normal fan of this polytopal  realization is a modification of the authors' \emph{velocity fan} and was found by OpenAI's Astra model.
\end{abstract}
\maketitle

\section{Categorical $n$-associahedra}
In \cite{b:2-associahedra}, the second author defined 2-associahedra, a family of posets arising from the study of functoriality for Fukaya categories, and conjectured that each one is the face poset of a convex polytope.  In \cite{backman2024higher}, the authors of the present note produced a complete fan realization of all $2$-associahedra called the \emph{velocity fan}.  Furthermore, we introduced categorical $n$-associahedra as a natural generalization of 2-associahedra and demonstrated that the velocity fan construction extends to that setting.

In the same work, we verified that these fans were polytopal, i.e. projective, for a special class of $n$-associahedra and stated that projectivity would be further investigated in future work.  After we posted a preprint of our work to arXiv, Vincent Pilaud coded up the $n$-associahedra and informed us, much to our surprise, that the velocity fan of the 2-associahedron $W_{1,1,1}$ is not polytopal.  

%We had a natural modification of the velocity fan, considered at the beginning of the project which we expected would fix the issue.  However, we recently found that the modified velocity fan also fails to be projective for $W_{1,1,1,1}$.   

In the present announcement we describe a modification of the velocity fan which is always polytopal. This modification was found by OpenAI's Astra model, and a full proof of polytopality will be presented in future work.  To keep this note short and accessible, we will focus on 2-associahedra, although the fan extends to all $n$-associahedra.

We begin with an informal definition of 2-associahedra.

\begin{definition}
Let $\mathcal{L}= \{\ell_1, \ldots, \ell_k\}$ be a nonempty collection of distinct vertical lines in $\mathbb{R}^2$ ordered from left to right.  Let $\mathcal{P} = \{p_{i,j}: 1\leq i \leq k , 1\leq j \leq m_i\}$ be a nonempty finite collection of distinct points such that the points $p_{i,j}$ lie on line $\ell_i$, and the points on a fixed line are ordered from bottom to top.  Suppose further that  $|\mathcal{L}|+|\mathcal{P}| \geq 3$. We refer to $\mathcal{X} = (\mathcal{L},\mathcal{P})$ as an \emph{arrangement of lines and points}.  
\end{definition}

We will define $\widetilde{\mathcal{X}}$ to be the collection of all arrangements obtainable from $\mathcal{X}$ by allowing individual lines and points to move continuously as follows:   lines are allowed to move horizontally (carrying their points with them) without meeting other lines, and points are allowed to move vertically on lines without meeting other points.  

We will describe the 2-associahedron as a combinatorial poset which encodes the boundary strata of a certain compactification of $\widetilde{\mathcal{X}}$ modulo translations and positive dilations \cite{b:witch-curves}.  Our approach is to first combinatorially describe the codimension-one strata, which we call \emph{collisions}, and then use a combinatorial rule to define the higher-codimension strata.  The collisions can be considered as a natural generalization of the facets of a classical associahedron, each of which corresponds to a consecutive set of points on a single line colliding and can equivalently be encoded by a pair of parentheses in a word (see \cite[Figure~1]{backman2024higher}).

\begin{definition}
A \emph{collision} for $\widetilde{\mathcal{X}}$ is a combinatorial object, which comes in two types, and encodes either
    \begin{enumerate}[itemsep=0pt, parsep=0pt, topsep=0pt]
    \item (type A): a consecutive collection of points on a single line coming together, or
    \item (type B): a consecutive collection of lines coming together, and the points on these lines coming together in groups.\footnote{We do not allow a collision where all points and lines come together simultaneously.}
\end{enumerate} 
\end{definition}

Forgetting the actual positions of lines and points, a type $A$ collision is encoded by a consecutive set of points on a line, while a combinatorial encoding of a type $B$ collision requires more information; this is an ordered partition of the points on the lines which respects the order of the points on individual lines.  We refer to each block of points in the partition associated to a type $B$ collision, together with the participating lines, as a \emph{2-bracket}.

\begin{figure}[ht]
\centering
% FIGURE SCALE: collision types.
\scalebox{1.6}{%
\includegraphics[width=5.59cm]{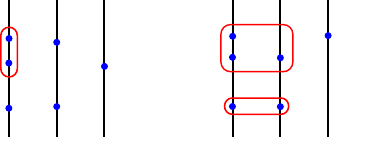}%
}
\caption{
\label{fig:collision_types}
A type $A$ collision (left) and a type $B$ collision (right). Here $\mathcal{X} = (\mathcal{L},\mathcal{P})$, with $\mathcal{L} = \{\ell_1, \ell_2, \ell_3\}$ and $\mathcal{P} = \{p_{1,1},p_{1,2},p_{1,3},p_{2,1},p_{2,2},p_{3,1}\}$.
}
\end{figure}

\begin{definition}
    We say that two collisions $\sC_1$ and $\sC_2$ are \emph{compatible} if
    
    \begin{enumerate}[itemsep=0pt, parsep=0pt, topsep=0pt]
        \item The sets of lines colliding in $\sC_1$ and $\sC_2$ are disjoint or nested.
        \item Any pair of 2-brackets in $\sC_1$ and $\sC_2$ are disjoint or nested. 
        \item The height partial orders on the points in $\sC_1$ and $\sC_2$ are compatible.
    \end{enumerate} 
\end{definition}

\begin{definition}
    A \emph{2-bracketing} for $\widetilde{\mathcal{X}}$ is an object which can be obtained as the union of a collection of compatible collisions.
\end{definition}

\begin{figure}[H]
\centering
% FIGURE SCALE: bracketing example.
\scalebox{1.2}{%
\includegraphics[width=4.19cm]{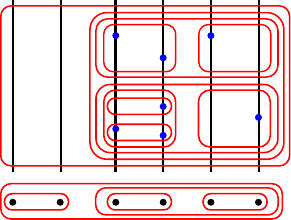}%
}
\caption{
\label{fig:2_bracketing}
A 2-bracketing.
}
\end{figure}

\begin{comment}
\begin{figure}[ht]
\centering
% FIGURE SCALE: nonexample (currently disabled).
\scalebox{1.1}{%
\includegraphics[width=4.19cm]{2-bracketing_non-examples.pdf}%
}
\caption{
\label{fig:2_bracketing_nonexample}
A diagram which is not a 2-bracketing because two groups of points intersect, but are not nested.
}
\end{figure}
\end{comment}

\begin{definition}
The \emph{2-associahedron} determined by $\widetilde{\mathcal{X}}$, which we denote $\cK(\widetilde{\mathcal{X}})$, is the collection of all 2-bracketings for $\widetilde{\mathcal{X}}$ ordered by refinement.  
\end{definition}

In keeping with the notation of \cite{b:2-associahedra}, we may also denote a 2-associahedron $\mathcal K(\widetilde{\mathcal X})$ by $W_{m_1,\ldots,m_k}$.

\begin{figure}[H]
\centering
% FIGURE SCALE: abstract W_{2,1}.
\scalebox{1.1}{%
\includegraphics[width=5.38cm]{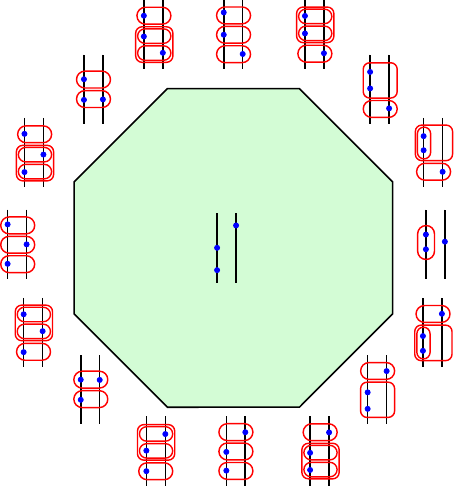}%
}
\caption{The $2$-associahedron $W_{2,1}$, with its faces labeled by $2$-bracketings.}
\label{fig:W_21_abstract}
\end{figure}

For an example of a 2-associahedron, see Figure~\ref{fig:W_21_abstract} -- the face poset of this octagon is a 2-associahedron.  It can occur that a 2-bracketing can be expressed as a union of collisions in multiple ways; this implies that 2-associahedra, unlike classical associahedra, are not always simple.

In \cite{backman2024higher}, we extend the above definition to arrangements of affine coordinate spaces in $\mathbb{R}^n$.  There, the combinatorial type of the arrangement is encoded by a rooted plane tree $\cT$ of depth $n$ (see \cite[Figure~3]{backman2024higher}), and we formally define collisions via natural maps of such trees.  Thus, we may denote an $n$-associahedron as $\cK(\cT)$ rather than $\cK(\widetilde{\mathcal{X}})$.

\section{Generalized velocity fans}
\label{velocitysection}

In this section, we recall the method for constructing a fan realization of 2-associahedra which was pursued in \cite{backman2024higher}, and we explain how to tweak that construction to produce a polytopal fan.

Let $\mathcal{X} = (\mathcal{L},\mathcal{P})$ be an arrangement of lines and points.  As noted above, the incidence relations for $\mathcal{X}$ are encoded by a rooted plane tree $\cT$.  As we are restricting to the case of $n=2$, the tree $\cT$ will have depth 2.  We begin by picking out a distinguished configuration inside of $\widetilde{\mathcal{X}}$ which we will denote $\mathcal{X}(\cT)$.  The key idea for the velocity fan is to interpret a collision $\sC$ as a transformation taking our distinguished arrangement $\mathcal{X}(\cT)$ to another distinguished arrangement $\mathcal{X}(\cT/\sC)$ over a single unit of time, and to encode the changes in relative distances between successive lines and points as entries in a vector $\rho({\sC})$, which will be our ray generator.  Here $\mathcal{X}(\cT/\sC)$ is the distinguished arrangement of points and lines where lines and points have come together according to the collision $\sC$.

 %Lines and points start from their positions in the distinguished configuration for $\widetilde{\mathcal{X}}$ and then move to another distinguished configuration for a smaller arrangement of lines and points $\widetilde{\mathcal{X}}/\sC$, where some of the lines and points have come together as described by $\sC$.

%The vector $\rho({\sC})$ will have $|\mathcal{L}|+|\mathcal{P}|-2$ entries.  The first $|\mathcal{L}|-1$ entries correspond to successive lines, and the second $|\mathcal{P}|-1$ entries correspond to successive points (with respect to the lexicographic order).  

\begin{definition}
Let $\sC$ be a collision.  Relabel the points: let $p_k$ be the $k$th point with respect to the lexicographic order on the pairs $(i,j)$.
We define a vector ${\rho}(\sC) \in \mathbb{R}^{m}$ where $m = |\mathcal{L}|+|\mathcal{P}| - 2$.  Denote
\begin{align}
{\rho}(\sC) 
=
(x_1,\dots, x_l,  y_1,\dots,y_p),
\end{align}
where $l = |\mathcal{L}|-1$ and $p = |\mathcal{P}|-1$.

We define the entries in ${\rho}(\sC)$ as follows:

\begin{equation}
\begin{aligned}
x_i = {\rm \,\,change\,\, in\,\, horizontal \,\, distance\,\, between\,\,} \ell_i {\rm \,\,and\,\,} \ell_{i+1}\\
{\rm\,\, when \,\,passing\,\, from\,\,} \mathcal{X}(\cT) {\rm \,\,to\,\,} \mathcal{X}(\cT/\sC),
\end{aligned}
\end{equation}

\begin{equation}
\begin{aligned}
y_i = {\rm \,\,change\,\, in\,\, vertical\,\, distance\,\, between\,\,} p_i {\rm \,\,and\,\,} p_{i+1}\\
    {\rm\,\, when \,\,passing\,\, from\,\,} \mathcal{X}(\cT) {\rm \,\,to\,\,} \mathcal{X}(\cT/\sC).\footnotemark
    \end{aligned}
\end{equation}
\footnotetext{All changes are measured as initial distances minus final distances.}

We define the corresponding ray to be
$${\tau}(\sC)
\coloneqq
\{ \lambda {\rho}(\sC)+ \gamma {\bf 1}: \lambda, \gamma \in \mathbb{R}, \lambda \geq 0\},$$
where ${\bf 1}$ is the all-ones vector.
Given a 2-bracketing $\sB \in \cK(\cT)$, we define the corresponding cone
\begin{align}
\tau(\sB)
\coloneqq
conv\{\tau(\sC): \sC {\rm \,\, a \,\,collision}, \, \sC \leq \sB\}.
\end{align}

\end{definition}

\begin{definition}
Given a rooted plane tree $\cT$ of depth $2$, we define the \emph{velocity fan} to be the collection of polyhedral cones
\begin{align}
\cF(\cT) = \{\tau(\sB):\sB \in \cK(\cT)\}
\end{align}
with lineality space $\langle 1 \rangle_{\mathbb{R}}$.
%\null\hfill$\triangle$
\end{definition}

\noindent

 {\bf The original velocity fan:}  For describing this fan, we simply need to define $\mathcal{X}(\cT)$ for each rooted plane tree $\cT$; since $\cT/\sC$ is also a rooted tree, this choice will also determine $\mathcal{X}(\cT/\sC)$.  We place line $\ell_i$ at position $x=i$, and we place point $p_{i,j}=p_k$ at position $(x,y) =(i,k)$.  Essentially, all the lines in $\mathcal{X}(\cT)$ get placed at distances one unit apart, and all of the points in $\mathcal{X}(\cT)$ get placed at heights one unit apart going from left to right. One can observe that the $x_i$ entries of ${\rho}(\sC)$ are always 0 or 1, while the $y_i$ entries of ${\rho}(\sC)$ can take other integer values.

\

We now describe the polytopal modification of the definition of the velocity fan.

\

 {\bf The barycentric velocity fan:}   Here we take the same choice of $\mathcal{X}(\cT)$.  However, we choose a terminal arrangement for each collision $\sC$ which is different from $\mathcal{X}(\cT/\sC)$.  For the terminal arrangement, we first permute the heights of the points in the participating lines according to an arbitrary linear extension of the height partial order associated to the collision\footnote{This rearrangement of points is a well-explored technique in \cite{backman2024higher} where it is referred to as a \emph{compatible $\cT$-shuffle for $\sC$}.}, and then contract the lines and points to the width and height \emph{barycenters} of the associated 1-brackets and 2-brackets, respectively.

The construction extends to arbitrary $n$ by applying the same reordering and barycentric terminal arrangement at each depth of the rooted plane tree. We announce the following result.

\begin{theorem}\label{thm:main}
Let $\cT$ be a rooted plane tree of depth $n$ and let $\cK(\cT)$ be the associated categorical $n$-associahedron. The barycentric velocity fan of $\cT$ is a projective fan whose face poset is $\cK(\cT)$.

\end{theorem}

We provide here a specific polytopal support function for the barycentric velocity fan which works for all $2$-associahedra and naturally generalizes to all $n$-associahedra.  Let $$S(m)=\sum_{j=0}^{m-1}j^2=\frac{m(m-1)(2m-1)}6.$$
For a collision $\sC$ involving an interval of $t$ lines and
2-bracket point sets $B_1,\ldots,B_q$, assign its ray generator for the barycentric velocity fan $\rho({\sC})$ the height
$$
h_\Lambda(\rho({\sC}))
=-\Lambda S(t)-t\sum_{j=1}^q S(|B_j|).
$$
This function extends linearly on the cones of the barycentric velocity fan, and for $\Lambda$ large enough, $h_\Lambda$ becomes strictly convex.\footnote{We note that this function depends only on the sizes of the brackets in $\sC$.  For realizing $W_{a,b}$ with $a,b\geq1$ and $a+b\geq4$ via the original velocity fan, it is impossible to achieve a polytopal support function which depends only on this data.}

\begin{figure}[H]
\centering
% FIGURE SCALE: all three collision rows together.
\scalebox{.85}{%
\begin{minipage}{13.97cm}
\centering
\includegraphics[width=13.97cm]{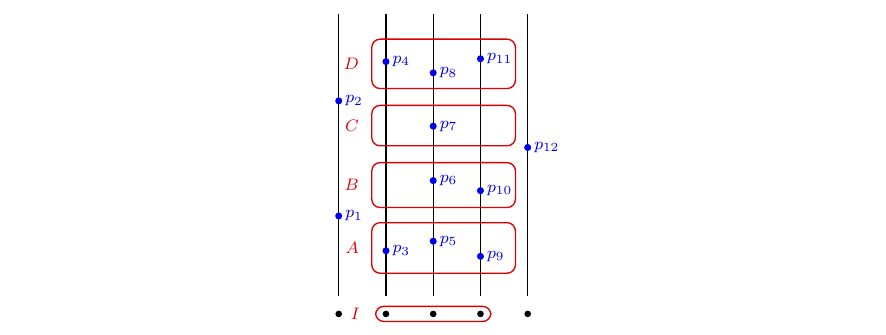}
\par\vspace{0.21cm} % Vertical gap after the abstract collision.
\includegraphics[width=13.97cm]{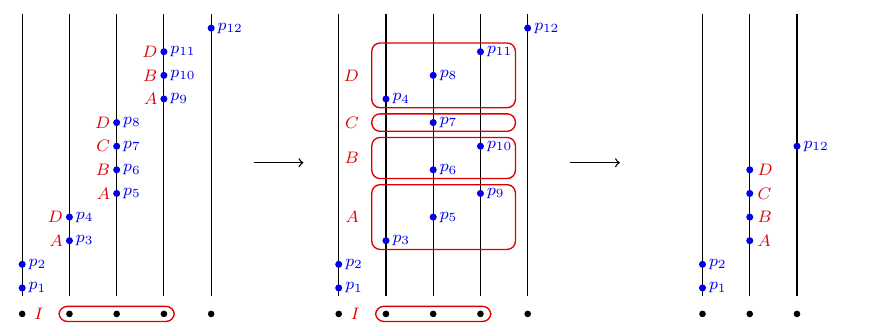}
\par\vspace{0.21cm} % Vertical gap between the two velocity rows.
\includegraphics[width=13.97cm]{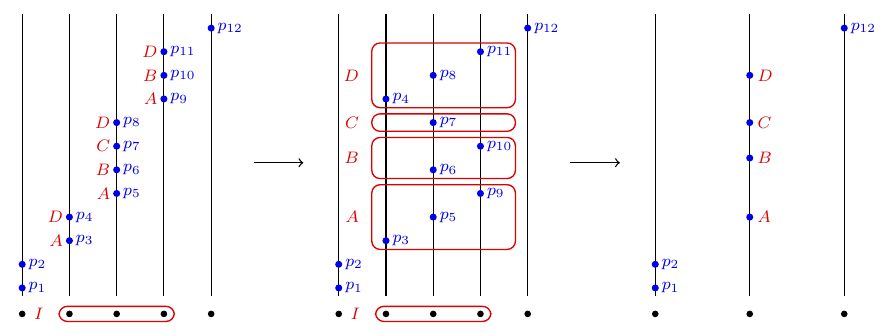}
\end{minipage}%
}
\caption{A collision in $W_{2,2,4,3,1}$ involving the middle three lines.
Top: the abstract collision.
Middle and bottom: the starting arrangement, the reordering of the participating points according to a linear extension of the height partial order, and the terminal arrangement for the original velocity fan and barycentric velocity fan, respectively.}
\label{fig:original_secondary_velocity_collision}
\end{figure}

For the collision $\sC$ in Figure~\ref{fig:original_secondary_velocity_collision},
we can compute the two ray generators explicitly.  The ray generator associated to the original velocity fan is 
\[
\rho_{\mathrm{orig}}(\sC)
=\bigl(0,1,1,0;\,0,0,-2,4,0,0,0,4,0,-1,0\bigr),
\] and the ray generator associated to the barycentric velocity fan is
\[
\rho_{\mathrm{bar}}(\sC)
=\bigl(-1,1,1,-1;\,0,-1,-5,7,-\tfrac32,-\tfrac12,-1,7,
        -\tfrac32,-\tfrac52,-1\bigr).
        \]
        
 \begin{figure}[H]
\centering
% FIGURE SCALE: both W_{2,1} fans together.
\scalebox{.95}{%
% The raisebox fractions align the fan origins.
\raisebox{\depth}{%
\raisebox{-0.45648189\height}{\includegraphics[width=6.32cm]{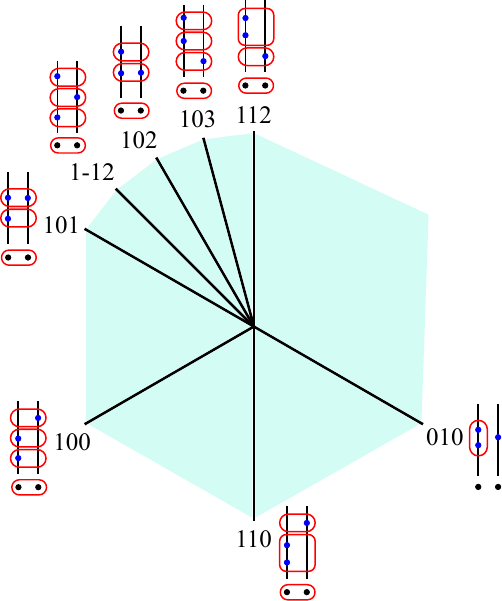}}%
\hspace{0.58cm}% Horizontal gap between the two fans.
\raisebox{-0.46195881\height}{\includegraphics[width=6.81cm]{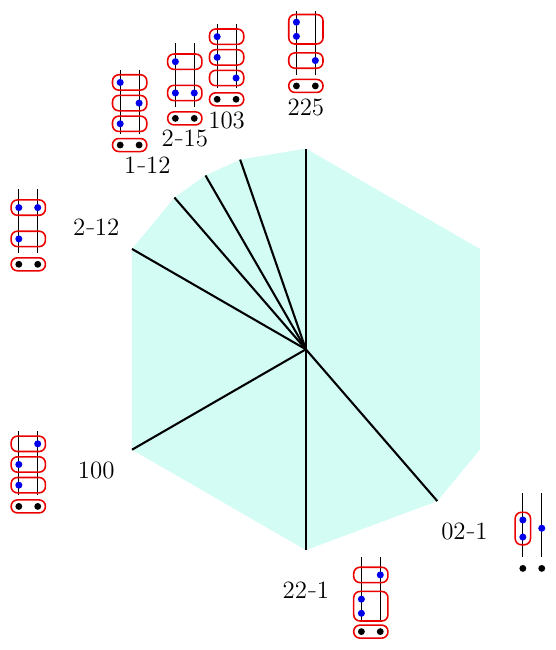}}%
}%
}
\caption{Two fan realizations of $W_{2,1}$: the original velocity fan (left)
and the barycentric velocity fan (right).}
\label{fig:W_21_velocity_fans_comparison}
\end{figure}

\begin{figure}[H]
\centering
% FIGURE SCALE: W_{2,2} net and both polytopes together.
\scalebox{.95}{%
% The raisebox commands keep all three components vertically centered.
\raisebox{\depth}{%
\raisebox{-0.5\height}{\includegraphics[width=6cm]{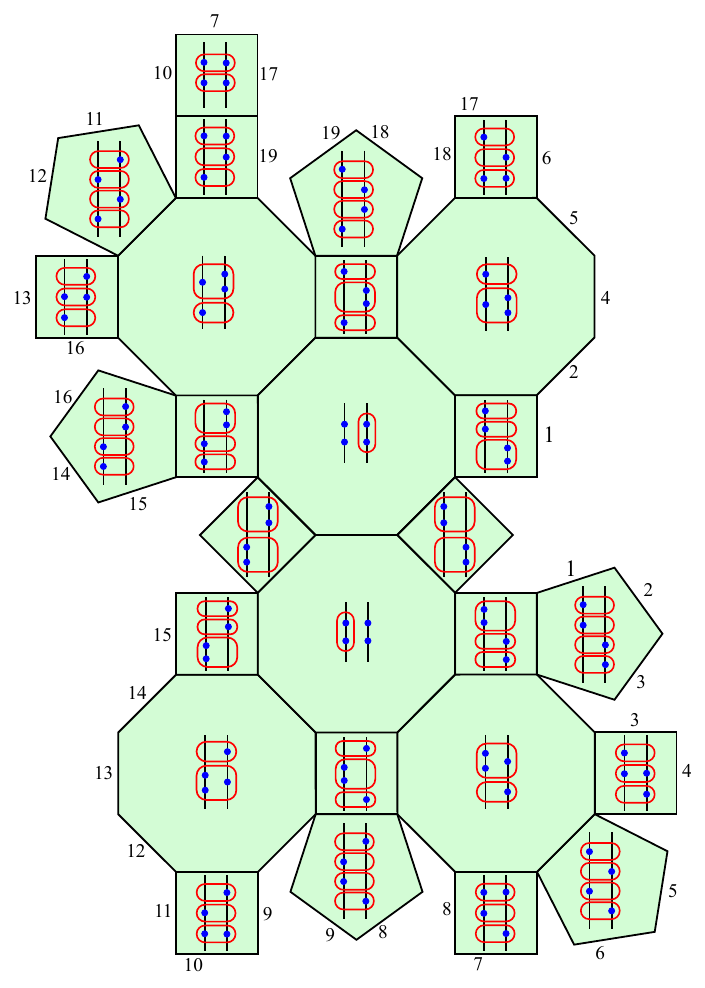}}%
\hspace{0.6cm}% Horizontal gap between the net and the original-fan polytope.
\raisebox{-0.5\height}{\includegraphics[height=8cm]{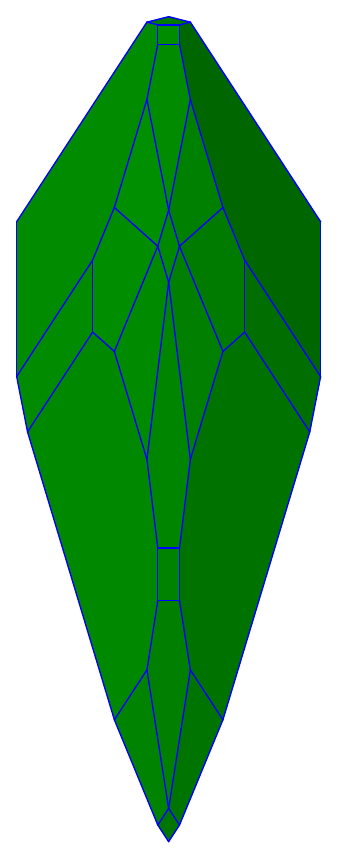}}%
\hspace{0.6cm}% Horizontal gap between the two polytopes.
\raisebox{-0.5\height}{\includegraphics[height=8cm]{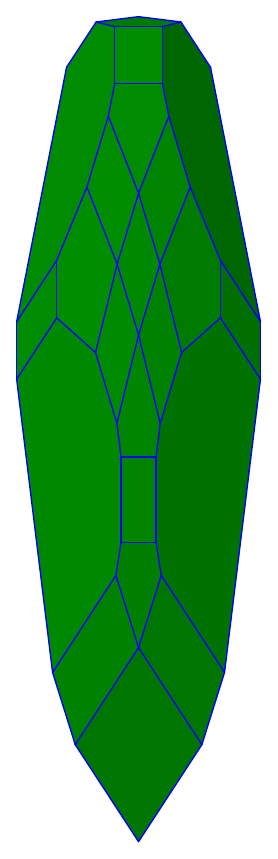}}%
}%
}
\caption{The $2$-associahedron $W_{2,2}$.
Left: a net representation.
Center: a polytopal realization with normal fan given by the original
velocity fan.
Right: a polytopal realization with normal fan given by the barycentric velocity fan.}
\label{fig:W22_net}
\label{fig:W22_original_and_secondary_realizations}
\end{figure}

\begin{figure}[H]
\centering
% FIGURE SCALE: W_{1,1,1} net and polytope together.
\scalebox{.8}{%
% Equal base heights keep the net and polytope balanced at any viewing angle.
\raisebox{\depth}{%
\raisebox{-0.5\height}{\includegraphics[height=7.42cm]{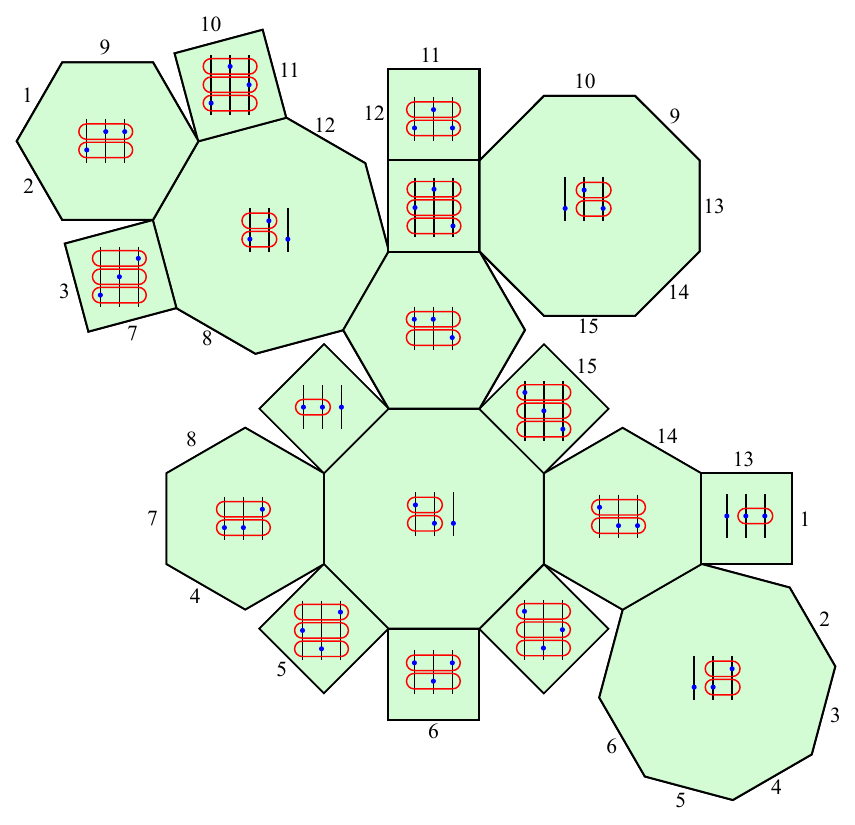}}%
\hspace{1cm}% Horizontal gap between the net and the polytope.
\raisebox{-0.5\height}{\includegraphics[height=7.42cm]{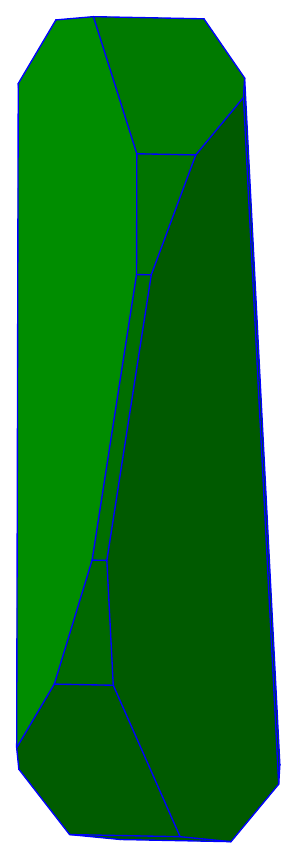}}%
}%
}
\caption{The $2$-associahedron $W_{1,1,1}$.
Left: a net representation.
Right: a polytopal realization whose normal fan is the barycentric velocity fan.  This is the only 3-dimensional 2-associahedron for which the original velocity fan is not polytopal.}
\label{fig:W111_secondary_realization}
\end{figure}

\section{Remarks}
\begin{itemize}
\item It is possible to upgrade the barycentric velocity fan by equipping the lines and points with distinct masses.
\item If the starting arrangement has only a single line, the velocity fan specializes to the normal fan of Loday's associahedron.  More generally, if the arrangement is \emph{concentrated}, meaning that all points lie on a single line, then the original velocity fan is a coarsening of the braid arrangement and is the normal fan of a generalized permutahedron.  

In contrast, the barycentric velocity fan for an arrangement with a single line recovers the secondary fan of a polygon with vertices on a parabola.  
%This interpretation of the fan was not originally communicated by Astra and was only found through further probing of the model.  
The secondary polytope of a polygon on a parabola has been investigated previously \cite{rote2003expansive,dimakis2011kp,galashin2022higher}, although this precise physical interpretation of the ray generators as relative velocities of points on a line colliding at their barycenters does not appear to have been noted previously.  More generally, the barycentric velocity fan of a concentrated 2-associahedron is a secondary fan closely related to those considered by Kerr--Palcic \cite{kerr2023painted}.  
%Interestingly, we had previously tried to understand whether there might be a relationship between the work of  \cite{kerr2023painted} and 2-associahedra, but had not made much progress in this direction.

\item  The original velocity fan admits a canonical smooth flag triangulation on the same set of rays -- this triangulation is motivated by the theory of wonderful compactifications.  In the case of the barycentric velocity fan, this triangulation is no longer smooth, but it is polytopal.

    \item There is a linear interpolation between the original and barycentric velocity fans:
    \[
    \rho_t(\mathcal C)=(1-t)\rho_{\mathrm{orig}}(\mathcal C)
+t\rho_{\mathrm{bar}}(\mathcal C),\qquad 0\leq t\leq1
\] consisting of fan realizations of $n$-associahedra.
\end{itemize}
\section*{AI Disclosure}
The authors posed the problem of producing a polytopal realization of $2$-associahedra to Astra on September 6th (the first day the model was publicly available) and received a positive response on the first try.  On September 8th, we received an email from Guillaume Laplante-Anfossi informing us that he had also given this task to Astra, which found the exact same fan, but provided different support function values for a polytopal realization.  

We remark that the fan we present here is slightly different from the fan which was found originally in that the lines are now chosen to also use barycentric terminal coordinates rather than the terminal coordinates of the original velocity fan.  This choice was made to allow for a uniform construction for all $n$.

We also used Astra in this note for generating figures and for copy editing, but not for writing the note itself.  We used Anthropic's Opus and Fable for producing a Lean 4 formalization that the barycentric velocity fan is polytopal for all categorical $n$-associahedra.

\section*{Acknowledgements}
We thank Vincent Pilaud for first coding up the categorical $n$-associahedra and informing us that the velocity fan of $W_{1,1,1}$ is not projective.  We thank Guillaume Laplante-Anfossi for informing us of his independent Astra run.  

S.B. was supported by a Simons Gift \#854037 and an NSF Grant
(DMS-2246967), and acknowledges support from OpenAI for the ChatGPT for Academic Researchers program.  D.P. was supported by the Deutsche Forschungsgemeinschaft (DFG, German Research
Foundation) – SFB-Geschäftszeichen 1624 – Projektnummer 506632645.

\bibliographystyle{alpha}
\bibliography{polytopal_realization}

\begin{thebibliography}{GPW22}

\bibitem[BBP24]{backman2024higher}
Spencer Backman, Nathaniel Bottman, and Daria Poliakova.
\newblock Higher-categorical associahedra.
\newblock arXiv preprint
  \href{https://arxiv.org/abs/2409.03633}{arXiv:2409.03633}, 2024.

\bibitem[Bot19a]{b:2-associahedra}
Nathaniel Bottman.
\newblock 2-associahedra.
\newblock {\em Algebraic \& Geometric Topology}, 19(2):743--806, 2019.

\bibitem[Bot19b]{b:witch-curves}
Nathaniel Bottman.
\newblock Moduli spaces of witch curves topologically realize the
  2-associahedra.
\newblock {\em Journal of Symplectic Geometry}, 17(6):1649--1682, 2019.

\bibitem[DMH11]{dimakis2011kp}
Aristophanes Dimakis and Folkert M{\"u}ller-Hoissen.
\newblock {KP} line solitons and {Tamari} lattices.
\newblock {\em Journal of Physics A: Mathematical and Theoretical},
  44(2):025203, 2011.

\bibitem[GPW22]{galashin2022higher}
Pavel Galashin, Alexander Postnikov, and Lauren Williams.
\newblock Higher secondary polytopes and regular plabic graphs.
\newblock {\em Advances in Mathematics}, 407:108549, 2022.

\bibitem[KP23]{kerr2023painted}
Gabriel Kerr and Sophia Palcic.
\newblock Painted tropical complexes.
\newblock arXiv preprint
  \href{https://arxiv.org/abs/2308.07409}{arXiv:2308.07409}, 2023.

\bibitem[RSS03]{rote2003expansive}
G{\"u}nter Rote, Francisco Santos, and Ileana Streinu.
\newblock Expansive motions and the polytope of pointed pseudo-triangulations.
\newblock In B.~Aronov, S.~Basu, J.~Pach, and M.~Sharir, editors, {\em Discrete
  and Computational Geometry: The {Goodman--Pollack} Festschrift}, volume~25 of
  {\em Algorithms and Combinatorics}, pages 699--736. Springer, Berlin,
  Heidelberg, 2003.

\end{thebibliography}

\end{document}